\documentclass[12pt]{article}

\usepackage[a4paper,margin=1in]{geometry}

\usepackage{amsmath,amsthm,amssymb,amsfonts,mathrsfs}
\usepackage{extarrows}
\usepackage{mathtools} 

\usepackage{tikz}
\usetikzlibrary{positioning}
\usetikzlibrary{calc}
\usetikzlibrary{arrows.meta,decorations.pathreplacing,calligraphy}

\usepackage{xcolor,enumerate,float}
\usepackage{booktabs}

\usepackage[numbers,square,comma,sort&compress]{natbib}
\usepackage[hidelinks]{hyperref}

\newtheorem{thm}{Theorem}[section]
\newtheorem{lem}[thm]{Lemma}
\newtheorem{cor}[thm]{Corollary}
\newtheorem{prop}[thm]{Proposition}

\theoremstyle{definition}
\newtheorem{defi}[thm]{Definition}

\newtheorem{rem}[thm]{Remark}

\newcommand{\defn}[1]{\emph{#1}}

\newcommand{\tre}[1]{\textcolor{red}{#1}}

\newcommand{\ga}{\gamma}

\newcommand{\sig}{\sigma}

\newcommand{\R}{\mathbb{R}}

\newcommand{\Z}{\mathbb{Z}}

\newcommand{\cD}{\mathcal{D}}

\newcommand{\cP}{\mathcal{P}}

\newcommand{\cU}{\mathcal{U}}
\newcommand{\cV}{\mathcal{V}}

\newcommand{\floor}[1]{\left\lfloor #1 \right\rfloor}

\newcommand{\what}{\widehat}

\newcommand{\BN}{\operatorname{BoNa}}
\newcommand{\diag}{\operatorname{diag}}

\newcommand{\minor}[2]{%
	\begin{pmatrix}
		#1 \\
		#2
	\end{pmatrix}	
}

\numberwithin{equation}{section}

\title{Analytic properties of Boolean-Narayana numbers}
\author{Xi Chen%
\thanks{
	Corresponding author.
	\newline\hspace*{3mm}
	{\it Email address:}\quad
	chenxi@dlut.edu.cn (X. Chen),
	sihao\_tang@hotmail.com (S. Tang)
},
Sihao Tang}
\date{\small
School of Mathematical Sciences, Dalian University of Technology, Dalian 116024, PR China}

\begin{document}
\maketitle

\begin{abstract}
Very recently, B\'ona introduced the Boolean-Narayana numbers $\BN(n,k)$ as the number of 0-1-trees on $n+1$ vertices with $k$ right edges.
We prove the total positivity of the Boolean-Narayana triangle by providing an explicit matrix factorization,
and construct a planar digraph for this triangle.
We give two combinatorial interpretations of the gamma-coefficients of the Boolean-Narayana polynomials.
We also show the strong $q$-log-convexity of the Boolean-Narayana polynomials and the asymptotic normality of the Boolean-Narayana numbers.
A bijection between 0-1-trees and colored Motzkin paths is constructed,
presenting an alternative interpretation of the Boolean-Narayana numbers.
\end{abstract}

\section{Introduction}

A \defn{binary plane tree} is a rooted plane tree on unlabeled vertices,
where every vertex has at most two children,
and each child is distinguished as a left child or a right child.
The $n$-th Catalan number $C_n=\binom{2n}{n}/(n+1)$ is known to count the number of binary plane trees on $n$ vertices.
The Narayana numbers $N(n,k)$ refine $C_{n+1}$ by counting binary plane trees on $n+1$ vertices with $k$ right edges.

If, in addition, every vertex that has two children is labeled $0$ or $1$,
then such a binary plane tree is called a \defn{0-1-tree}.
The number of 0-1-trees on $n+1$ vertices is called the \defn{$n$th Boolean-Catalan number} \cite{Bon26b}.
Very recently, B\'ona \cite{Bon26a} refined these numbers by recording the number of right edges.

\begin{defi}\label{defi:BN}
	The \defn{Boolean-Narayana numbers}, denoted by $\BN(n,k)$,
	are defined as the number of 0-1-trees on $n+1$ vertices with $k$ right edges.
\end{defi}

\begin{figure}
	\centering
	\begin{tikzpicture}[scale=0.7,font=\footnotesize,
		vtx/.style={circle,fill=black,inner sep=1.5pt},
		]
		\begin{scope}
			\node[vtx] (a) at (0,0) {};
			\node[vtx] (b) at (-1,-1) {};
			\node[vtx] (c) at (-2,-2) {};
			\draw (a)--(b)--(c);
		\end{scope}
		\draw[decorate,decoration={brace}] (0,-2.5)--(-2,-2.5) node[midway,below=2pt]{$\BN(2,0)=1$};
		
		\begin{scope}[xshift=2cm]
			\node[vtx] (a) at (0,0) {};
			\node[vtx] (b) at (-1,-1) {};
			\node[vtx] (c) at (0,-2) {};
			\draw (a)--(b)--(c);
		\end{scope}
		
		\begin{scope}[xshift=3.5cm]
			\node[vtx] (a) at (0,0) {};
			\node[vtx] (b) at (1,-1) {};
			\node[vtx] (c) at (0,-2) {};
			\draw (a)--(b)--(c);
		\end{scope}

		\begin{scope}[xshift=7cm]
			\node[vtx,label=above:{$0$}] (a) at (0,0) {};
			\node[vtx] (b) at (-1,-1) {};
			\node[vtx] (c) at (1,-1) {};
			\draw (b)--(a)--(c);
		\end{scope}

		\begin{scope}[xshift=10cm]
			\node[vtx,label=above:{$1$}] (a) at (0,0) {};
			\node[vtx] (b) at (-1,-1) {};
			\node[vtx] (c) at (1,-1) {};
			\draw (b)--(a)--(c);
		\end{scope}
		\draw[decorate,decoration={brace}] (11,-2.5)--(1,-2.5) node[midway,below=2pt]{$\BN(2,1)=4$};

		\begin{scope}[xshift=12cm]
			\node[vtx] (a) at (0,0) {};
			\node[vtx] (b) at (1,-1) {};
			\node[vtx] (c) at (2,-2) {};
			\draw (a)--(b)--(c);
		\end{scope}
		\draw[decorate,decoration={brace}] (14,-2.5)--(12,-2.5) node[midway,below=2pt]{$\BN(2,2)=1$};
	\end{tikzpicture}
	\caption{The 0-1-trees on three vertices.}
	\label{fig:0-1-tree}
\end{figure}
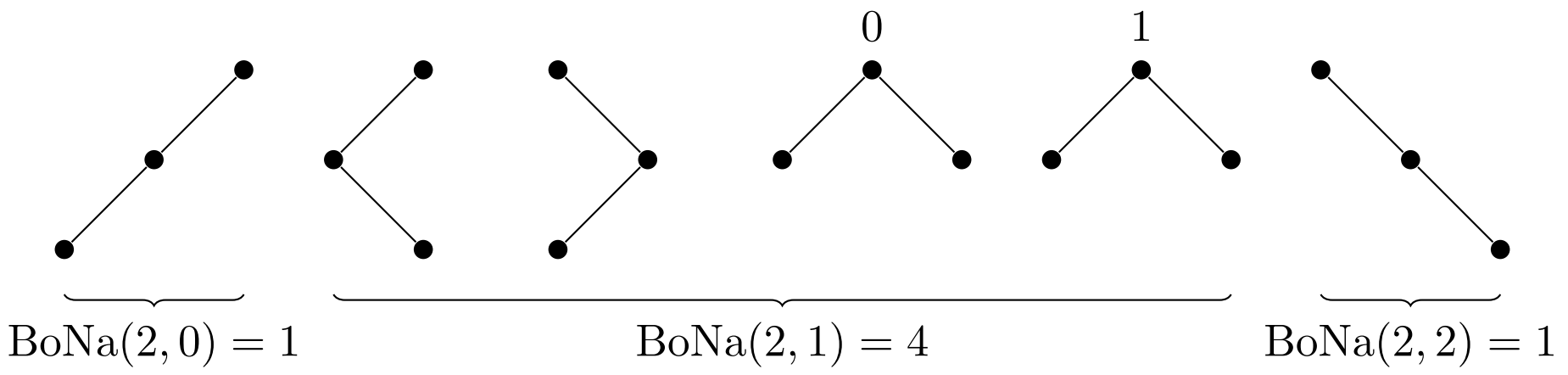

Note that our indices $n$ and $k$ start at zero,
while those in B\'ona's \cite{Bon26a} start at one.
Figure \ref{fig:0-1-tree} shows the 0-1-trees on three vertices.
B\'ona \cite{Bon26a} proved the symmetry $\BN(n,k)=\BN(n,n-k)$
and gave the explicit formula
\begin{equation}\label{eq:BN-formula}
	\BN(n,k)=\frac{1}{n+1} \binom{n+1}{k} \sum_{j=0}^k 2^j \binom{k}{j} \binom{n-k+1}{j+1},
	\quad 0\le k\le n/2.
\end{equation}
We call the infinite lower triangular matrix
\[
\BN:=[\BN(n,k)]_{n,k\ge0}=
\begin{bmatrix}
	1 \\
	1 & 1 \\
	1 & 4 & 1 \\
	1 & 9 & 9 & 1 \\
	1 & 16 & 38 & 16 & 1 \\
	\vdots &&&&& \ddots
\end{bmatrix}
\]
the \defn{Boolean-Narayana triangle}.
The \defn{Boolean-Narayana polynomials}
\[
\BN_n(q)=\sum_{k=0}^n \BN(n,k)q^k
\]
are defined as the row-generating functions of the Boolean-Narayana triangle.
Very recently, B\'ona \cite{Bon26a} proved the real-rootedness and interlacing property of $\BN_n(q)$,
and the log-concavity of $\BN(n,k)$ for fixed $n$ or fixed $k$.
This paper is devoted to some other analytic properties of the Boolean-Narayana numbers.

In Section 2, we prove that the Boolean-Narayana triangle $\BN$ is totally positive by providing a totally positive matrix factorization,
and construct a planar digraph for $\BN$.
In Section 3, we interpret the $\ga$-coefficients of $\BN_n(q)$ in terms of weighted 0-1-trees.
Then, using the explicit $\ga$-expansion, we prove that $\BN_n(q)$ form a strongly $q$-log-convex sequence.
In Section 4, we obtain the asymptotic normality of the Boolean-Narayana numbers.
In Section 5, we construct a bijection between 0-1-trees and colored Motzkin paths,
hence present an alternative definition of the Boolean-Narayana numbers and an alternative interpretation of the $\ga$-coefficients of $\BN_n(q)$.

\section{Total positivity}

Following Karlin \cite{Kar68}, a finite or infinite matrix is called \defn{totally positive} (TP for short) if its minors of all orders are nonnegative;
it is called \defn{TP$_r$} if its minors of orders $\le r$ are nonnegative.
Totally positive matrices have been widely studied in both analysis and combinatorics \cite{Bre95,Kar68,Pin10}.
Over the past five years, total positivity of combinatorial matrices has attracted considerable attention \cite{Zhu21,LLY+22,Sok23,CS24,Zhu24,MWZ26,Zhu26}.
We review some basic facts about TP matrices.
The first proposition follows from the classic Cauchy-Binet formula,
and the second follows directly from the definition.

\begin{prop}\label{prop:TP-prod}
	The product of two TP matrices is still TP.
\end{prop}

\begin{prop}\label{prop:TP-submat}
	Every submatrix of a TP matrix is also TP.
\end{prop}

Many combinatorial triangles have been shown to be totally positive.
For example, it is well known that the Pascal triangle is TP (see \cite[p. 137]{Kar68} and \cite[Remark 2.4]{CLW15a} for details).
The first author et al. \cite{CWZ20} studied the total positivity of combinatorial triangles related to the Motzkin numbers.
Wang and Yang \cite{WY18} showed that the Narayana triangle is TP.
We obtain the following TP result for the Boolean-Narayana triangle.

\begin{thm}\label{thm:TP}
	The Boolean-Narayana triangle $\BN$ is totally positive.
\end{thm}

Our proof technique is to decompose $\BN$ into several TP matrices.
Then, by Proposition \ref{prop:TP-prod}, the total positivity of $\BN$ follows.
The proof of Theorem \ref{thm:TP} is divided into the following two lemmas.

\begin{lem}\label{lem:BN-factor}
	We have
	\[
	\BN=PD_1QD_2,
	\]
	where $P=\left[\binom{n}{k}\right]_{n,k\ge0}$ is the Pascal triangle,
	\[
	D_1=\diag(2^n n!)_{n\ge0},\quad
	D_2=\diag\left(\frac{1}{2^k k!(k+1)!}\right)_{k\ge0},
	\]
	and $Q=[Q_{n,k}]_{n,k\ge0}$ is the matrix defined by
	\[
	Q_{n,k}=
	\begin{cases}
		\dfrac{k!(k+1)!}{(2k-n)!(n-k)!(n-k+1)!}, & k\le n\le2k;\\
		0, &\text{otherwise}.
	\end{cases}
	\]
\end{lem}

\begin{proof}
By \eqref{eq:BN-formula} and the symmetry $\BN(n,k)=\BN(n,n-k)$, we have
\begin{align}
	\nonumber
	\BN(n,k)
	&=\frac{1}{n+1} \binom{n+1}{k} \sum_{j=0}^{\min\{k,n-k\}} 2^j \binom{k}{j} \binom{n-k+1}{j+1}\\
	\label{eq:BN-formula-2}
	&=\sum_{j=0}^{\min\{k,n-k\}} \binom{n}{k+j} \binom{k+j}{k-j} 2^j C_j,
\end{align}
where $C_j=\binom{2j}{j}/(j+1)$ is the $j$-th Catalan number.
Replacing $j$ by $m-k$ gives
\begin{align*}
	\BN(n,k)
	&=\sum_{m=k}^{\min\{2k,n\}} \binom{n}{m} \binom{m}{2k-m} 2^{m-k} C_{m-k}\\
	&=\sum_{m=k}^{\min\{2k,n\}} \binom{n}{m}\cdot 2^m m!\cdot \frac{k!(k+1)!}{(2k-m)!(m-k)!(m-k+1)!}\cdot \frac{1}{2^k k!(k+1)!}.
\end{align*}
The last equation implies that $\BN=PD_1QD_2$.
\end{proof}

\begin{rem}
Note that there is a common factor $k!(k+1)!$ in $Q$ and in the denominators of $D_2$.
We suggest not cancelling this factor so that a planar digraph with integer weights can be constructed
(Corollary \ref{cor:Q-network}).
\end{rem}

Recall that the Pascal triangle $P$ is TP.
Clearly, a diagonal matrix with nonnegative diagonal entries is TP.
So both $D_1$ and $D_2$ are TP.
It remains to prove the total positivity of $Q$.
We will embed $Q$ into a TP matrix of a larger size,
then Proposition \ref{prop:TP-submat} applies.
  
\begin{lem}\label{lem:L}
	Let $L=[L_{r,s}]_{r,s\ge0}$ be the matrix defined by
	\begin{equation}\label{eq:W}
		L_{r,s}=
		\begin{cases}
			Q_{n,k}, & r=2n, s=2k, k\le n\le2k;\\
			\frac{k+2}{n-k+1} Q_{n-1,k}, & r=2n, s=2k+1, k+1\le n\le2k+1;\\
			0, & \text{otherwise}.
		\end{cases}
	\end{equation}
	Then $L$ is totally positive.
	Furthermore, the matrix $Q=[Q_{n,k}]_{n,k\ge0}$ is also totally positive.
\end{lem}

\begin{proof}
By definition $Q_{n,k}=L_{2n,2k}$ for all $n,k\ge0$.
That is, $Q$ is the submatrix of $L$ obtained by taking all even-indexed rows and columns.
So by Proposition \ref{prop:TP-submat}, it suffices to prove that $L$ is TP.

For $i\ge1$, let $A_i$ and $C_i$ be the block matrices
\[
A_i:=
\begin{bmatrix}
	I_{2i-1} & & & \\
& \boxed{\begin{matrix}
		i-1 & 0\\
		1   & i+1
\end{matrix}} & & \\[4pt]
& & \boxed{\begin{matrix}
		i-1 & 0\\
		1   & i+1
\end{matrix}} & \\[4pt]
& & & \ddots
\end{bmatrix}
\]
and
\[
C_i:=
\begin{bmatrix}
I_{2i} & & & \\
& \boxed{\begin{matrix}
		1 & 0\\
		1 & 1
\end{matrix}} & & \\[4pt]
& & \boxed{\begin{matrix}
		1 & 0\\
		1 & 1
\end{matrix}} & \\[4pt]
& & & \ddots
\end{bmatrix},
\]
respectively.
Let $T_i=A_i C_i$.
Then
\begin{equation}\label{eq:Ti}
	T_i=
	\begin{bmatrix}
		I_{2i-1} \\
				& \what{T}_i
	\end{bmatrix},
\end{equation}
where $\what{T}_i$ is
lower triangular with two subdiagonals.
The entries on the three (sub)diagonals alternate between $i-1$ and $i+1$,
between $1$ and $i-1$, and between $0$ and $1$, respectively:
\[
\what{T}_i=
\begin{bmatrix}
	i-1 \\
	1   & i+1 \\
	0   & i-1 & i-1 \\
	    & 1   & 1   & i+1 \\
		&     & 0   & i-1 & i-1    \\
		&     &     & 1   & 1      & i+1 \\
		&     &     &     & \ddots & \ddots & \ddots
\end{bmatrix}.
\]
Clearly, the bidiagonal matrices $A_i$ and $C_i$ are TP.
Hence, by Proposition \ref{prop:TP-prod}, the matrices $T_i$ are still TP.

In what follows, we will prove that
\[
L=\prod_{i\ge1} T_i=\prod_{i\ge1} (A_i C_i).
\]
Let $M_0=I$ and $M_j=T_1 T_2\cdots T_j$ for $j\ge 1$.
By the definitions of $T_j$, for each fixed $s\ge0$,
the $s$th column of $T_j$ coincides with that of $I$
for all sufficiently large $j$.
Thus the $s$th column of $M_j=M_{j-1}T_j$ eventually
remains unchanged.
Consequently, every entry of $M_j$ is eventually constant,
and the limit
\[
M:=\lim_{j\to\infty}M_j
\]
exists entrywise.
Each finite minor of $M$ agrees with the corresponding minor of some $M_j$,
since the finitely many columns involved have stabilized.
By Proposition \ref{prop:TP-prod}, each $M_j$ is TP.
Hence, $M$ is also TP.
It remains to show that $M=L$.

Note that all the matrices we consider have row/column indices starting at zero.
The odd-indexed rows of $A_1$ vanish,
so the corresponding rows of $M=A_1C_1\cdots$ vanish as well.
Thus, $M_{2n+1,s}=0=L_{2n+1,s}$.
For the even-indexed rows of $M$,
since $M_j=M_{j-1}T_j$, we have
\begin{align*}
	(M_j)_{2n,2k+1}
	&=\sum_{\ell} (M_{j-1})_{2n,\ell} (T_j)_{\ell,2k+1}\\
	&=(M_{j-1})_{2n,2k+1} (T_j)_{2k+1,2k+1}
	+(M_{j-1})_{2n,2k+2} (T_j)_{2k+2,2k+1}\\
	&=(j-1)(M_{j-1})_{2n,2k+1}+(M_{j-1})_{2n,2k+2}
\end{align*}
for $n\ge j$ and $k\ge j-1$,
and
\begin{align*}
	(M_j)_{2n,2k}
	&=\sum_{\ell} (M_{j-1})_{2n,\ell} (T_j)_{\ell,2k}\\
	&=(M_{j-1})_{2n,2k} (T_j)_{2k,2k}
	+(M_{j-1})_{2n,2k+1} (T_j)_{2k+1,2k}
	+(M_{j-1})_{2n,2k+2} (T_j)_{2k+2,2k}\\
	&=(j+1)(M_{j-1})_{2n,2k}+(j-1)(M_{j-1})_{2n,2k+1}+(M_{j-1})_{2n,2k+2}
\end{align*}
for $n\ge j$ and $k\ge j$.
Induction on $j$ gives, for the odd-indexed columns with $k\ge j-1$,
\[
(M_j)_{2n,2k+1}=
\begin{cases}
	\dfrac{(j-1)!(j+1)!}{(n-k-1)!(n-k+1)!(j-n+k)!}, & k+1\le n\le k+j;\\
	0, & \text{otherwise}
\end{cases}
\]
and, for the even-indexed columns with $k\ge j$,
\[
(M_j)_{2n,2k}=
\begin{cases}
	\dfrac{j!(j+1)!}{(n-k)!(n-k+1)!(j-n+k)!}, & k\le n\le k+j;\\
	0, & \text{otherwise}.
\end{cases}
\]
By \eqref{eq:Ti}, columns $2k+1$ and $2k$ stabilize at $M_{k+1}$ and $M_k$, respectively, with $M_0=I$.
Hence,
\[
M_{2n,2k+1}=(M_{k+1})_{2n,2k+1}
=\frac{k!(k+2)!}{(2k+1-n)!(n-k-1)!(n-k+1)!}=L_{2n,2k+1}
\]
for $k+1\le n\le 2k+1$, and
\[
M_{2n,2k}=(M_k)_{2n,2k}
=\frac{k!(k+1)!}{(2k-n)!(n-k)!(n-k+1)!}=L_{2n,2k}
\]
for $k\le n\le 2k$,
and they both equal $0$ otherwise.
Therefore, $M=L$, as desired.
\end{proof}

The proof of Lemma \ref{lem:L} also gives rise to a planar digraph characterization of the matrices $L$ and $Q$.
Then the famous Lindstr\"om-Gessel-Viennot (LGV) lemma \cite[Chapter 32]{AZ18} also implies their total positivity.

Let $\cD=(V,E)$ be an acyclic planar digraph,
in which there are finitely many paths between any two vertices.
Each edge is equipped with a weight.
The weight $w(\cP)$ of a path $\cP$ is defined as the product of the weights of the edges traversed by $\cP$.
For $u, v\in V$, define
$$
P_{\cD}(u\to v):=\sum_{\cP:\,u\to v} w(\cP),
$$
where the sum ranges over all directed paths going from $u$ to $v$ in $\cD$.
We adopt the convention that $P_{\cD}(u\to u)=1$.
Let $\cU=\{u_0,u_1,u_2,\ldots\}$ and $\cV=\{v_0,v_1,v_2,\ldots\}$
be the sets of \defn{sources} and \defn{sinks} of $\cD$, respectively.
Define the \defn{path matrix} of $\cD$ as
\[
P_{\cD}:=\left[P_{\cD}(u_n\to v_k)\right]_{n,k\ge0}.
\]
The LGV lemma shows that every minor of $P_{\cD}$ is a sum over families of nonintersecting paths between specific subsets of $\cU$ and $\cV$:
$$
P_{\cD}\minor{n_1,\ldots,n_r}{k_1,\ldots,k_r}=
\sum_{\substack{\sigma\in\mathfrak{S}_r \\
\cP_i:\,u_{n_i}\to v_{k_{\sigma(i)}} \\
\cP_1,\ldots,\cP_r ~\text{nonintersecting}}}
{\rm sgn}(\sigma) w(\cP_1)\cdots w(\cP_r).
$$
It follows that if the edge weights of $\cD$ are all nonnegative,
and there are no nonintersecting paths unless $\sig$ is the identity,
then the path matrix $P_{\cD}$ is totally positive.

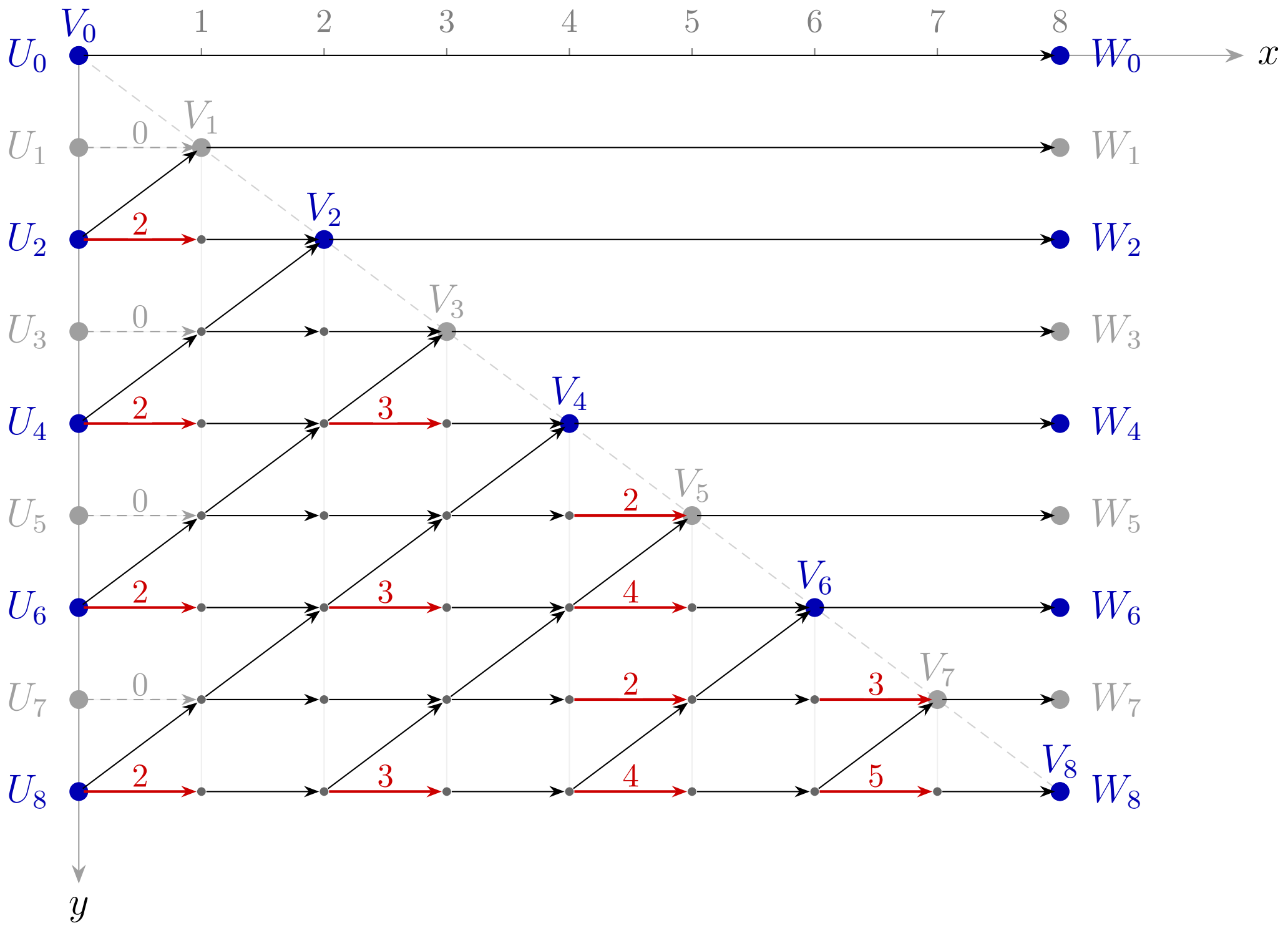
\begin{figure}[htbp]
	\centering
	\begin{tikzpicture}[x=1cm,y=-0.75cm,scale=1.3,
		axis/.style={-{Stealth[length=2mm]},thin,gray!75},
		edge/.style={-{Stealth[length=1.6mm]},shorten <=1.5pt,shorten >=1.5pt},
		nonunit/.style={edge,red!80!black,thick},
		zerowt/.style={edge,gray!75},
		weight/.style={fill=white,inner sep=1.1pt,font=\footnotesize}
		]
		\draw[axis] (0,0)--(9.5,0) node[right,black] {$x$};
		\draw[axis] (0,0)--(0,9) node[below,black] {$y$};
		\foreach \r in {1,...,8}{
			\draw[gray] (\r,0)--++(0,-.08);
			\node[gray,above=3pt,font=\footnotesize] at (\r,0) {$\r$};
			\draw[gray!12,thin] (0,\r)--(\r,\r);
			\draw[gray!12,thin] (\r,\r)--(\r,8);
		}
		\draw[gray!35,densely dashed] (0,0)--(8,8);
		\foreach \y in {0,...,8}{
			\foreach \x in {0,...,\y}{
				\fill[black!60] (\x,\y) circle[radius=1pt];
			}
		}
		\foreach \r in {0,...,8}{
			\pgfmathtruncatemacro{\rparity}{mod(\r,2)}
			\ifnum\rparity=0\relax
				\def\sourcecolor{blue!70!black}
				\def\sinkcolor{blue!70!black}
			\else
				\def\sourcecolor{gray!75}
				\def\sinkcolor{gray!75}
			\fi
			\node[circle,fill=\sourcecolor,inner sep=2pt] at (0,\r) {};
			\node[left=5pt,text=\sourcecolor] at (0,\r) {$U_{\r}$};
			\node[circle,fill=\sinkcolor,inner sep=2pt] at (\r,\r) {};
			\node[above,text=\sinkcolor] at (\r,\r) {$V_{\r}$};
			\node[circle,fill=\sinkcolor,inner sep=2pt] at (8,\r) {};
			\node[right=5pt,text=\sinkcolor] at (8,\r) {$W_{\r}$};	
		}
		\foreach \x in {0,...,7}{
			\foreach \y in {0,...,8}{
				\ifnum\y>\x\relax
					\pgfmathtruncatemacro{\xp}{\x+1}
					\pgfmathtruncatemacro{\xparity}{mod(\x,2)}
					\pgfmathtruncatemacro{\yparity}{mod(\y,2)}
					\ifnum\xparity=1\relax
						\def\edgeweight{1}
					\else
						\ifnum\yparity=0\relax
							\pgfmathtruncatemacro{\edgeweight}{\x/2+2}
						\else
							\pgfmathtruncatemacro{\edgeweight}{\x/2}
						\fi
					\fi
					\ifnum\edgeweight>0\relax
						\ifnum\edgeweight=1\relax
							\draw[edge] (\x,\y)--(\xp,\y);
						\else
							\draw[nonunit] (\x,\y)--node[weight,above] {$\edgeweight$}(\xp,\y);
						\fi
					\else
						\draw[zerowt,densely dashed] (\x,\y)--node[weight,above] {$\edgeweight$}(\xp,\y);
					\fi
					\ifnum\xparity=\yparity\relax
						\pgfmathtruncatemacro{\ym}{\y-1}
						\ifnum\ym<\xp\relax
						\else
							\draw[edge] (\x,\y)--(\xp,\ym);
						\fi
					\fi
				\fi
			}
		}
		\foreach \r in {0,...,7}{
			\draw[edge] (\r,\r)--(8,\r);	
		}
	\end{tikzpicture}
	\caption{The weighted digraph $\cD$ restricted to $0\le x\le y\le8$.
	All unlabeled edges have weight $1$.
	Even-indexed sources and sinks give the leading $5\times5$ submatrix of $Q$.}
	\label{fig:network-Q}
\end{figure}

\begin{cor}\label{cor:Q-network}
Let $\cD$ be a weighted digraph on the vertex set
\[
\{(x,y)\in\Z^2: y\ge x\ge0\}
\]
with the following edges whenever both endpoints lie in the vertex set:
\begin{enumerate}[\rm(i)]
	\item if $x$ and $y$ are both odd,
	put a horizontal edge $(x,y)\to(x+1,y)$ and an upward diagonal edge $(x,y)\to(x+1,y-1)$, 
	both with weight $1$;
	\item if $x$ is odd and $y$ is even,
	put a horizontal edge $(x,y)\to(x+1,y)$ with weight $1$;
	\item if $x$ is even and $y$ is odd, 
	put a horizontal edge $(x,y)\to(x+1,y)$ with weight $x/2$;
	\item if $x$ and $y$ are both even, 
	put a horizontal edge $(x,y)\to(x+1,y)$ and an upward diagonal edge $(x,y)\to(x+1,y-1)$, 
	with weights $x/2+2$ and $1$,
	respectively.
\end{enumerate}
Then the path matrix of $\cD$ with sources
$U_{2n}=(0,2n)$ and sinks $V_{2k}=(2k,2k)$ is precisely the matrix $Q=[Q_{n,k}]_{n,k\ge0}$ defined in Lemma \ref{lem:BN-factor}.
(See Figure \ref{fig:network-Q} for instance.)
\end{cor}

\begin{figure}[htbp]
	\centering
	\begin{tikzpicture}[x=1cm,y=-0.75cm,scale=1.3,font=\small,
		vertex/.style={circle,fill=black,inner sep=1.5pt},
		edge/.style={-{Stealth[length=1.6mm]},shorten <=1.5pt,shorten >=1.5pt},
		nonunit/.style={edge,red!80!black,thick},
		weight/.style={fill=white,inner sep=1.1pt,font=\footnotesize}
		]
		\begin{scope}
			\foreach \y/\r in {0/0,1/{2i-2},2/{2i-1},3/{2i},4/{2i+1},5/{2i+2},6/{2i+3}}{
				\node[vertex,label=left:{$\r$}] at (0,\y) {};
				\node[vertex] at (1,\y) {};
				\draw[edge] (0,\y)--(1,\y);
			}
			\foreach \y in {3,5}{
				\draw[edge] (0,\y)--(1,\y-1);
				\draw[nonunit] (0,\y-1)--node[weight,above] {$i-1$}(1,\y-1);
				\draw[nonunit] (0,\y)--node[weight,above] {$i+1$}(1,\y);
			}
			\draw[nonunit] (0,6)--node[weight,above] {$i-1$}(1,6);
			\foreach \y in {0.5,6.5}
				\foreach \x in {0,1}
					\node at (\x,\y) {$\vdots$};
			\node at (0.5,-0.5) {$A_i$};
		\end{scope}
		\begin{scope}[xshift=4cm]
			\foreach \y/\r in {0/0,1/{2i-2},2/{2i-1},3/{2i},4/{2i+1},5/{2i+2},6/{2i+3}}{
				\node[vertex,label=left:{$\r$}] at (0,\y) {};
				\node[vertex] at (1,\y) {};
				\draw[edge] (0,\y)--(1,\y);
			}
			\foreach \y in {4,6}{
				\draw[edge] (0,\y)--(1,\y-1);
			}
			\foreach \y in {0.5,6.5}
				\foreach \x in {0,1}
					\node at (\x,\y) {$\vdots$};
			\node at (0.5,-0.5) {$C_i$};
		\end{scope}
	\end{tikzpicture}
	\caption{Digraphs corresponding to $A_i$ and $C_i$.
	All unlabeled edges have weight $1$.
	Labels of the vertices indicate their ordinates.}
	\label{fig:network-AC}
\end{figure}
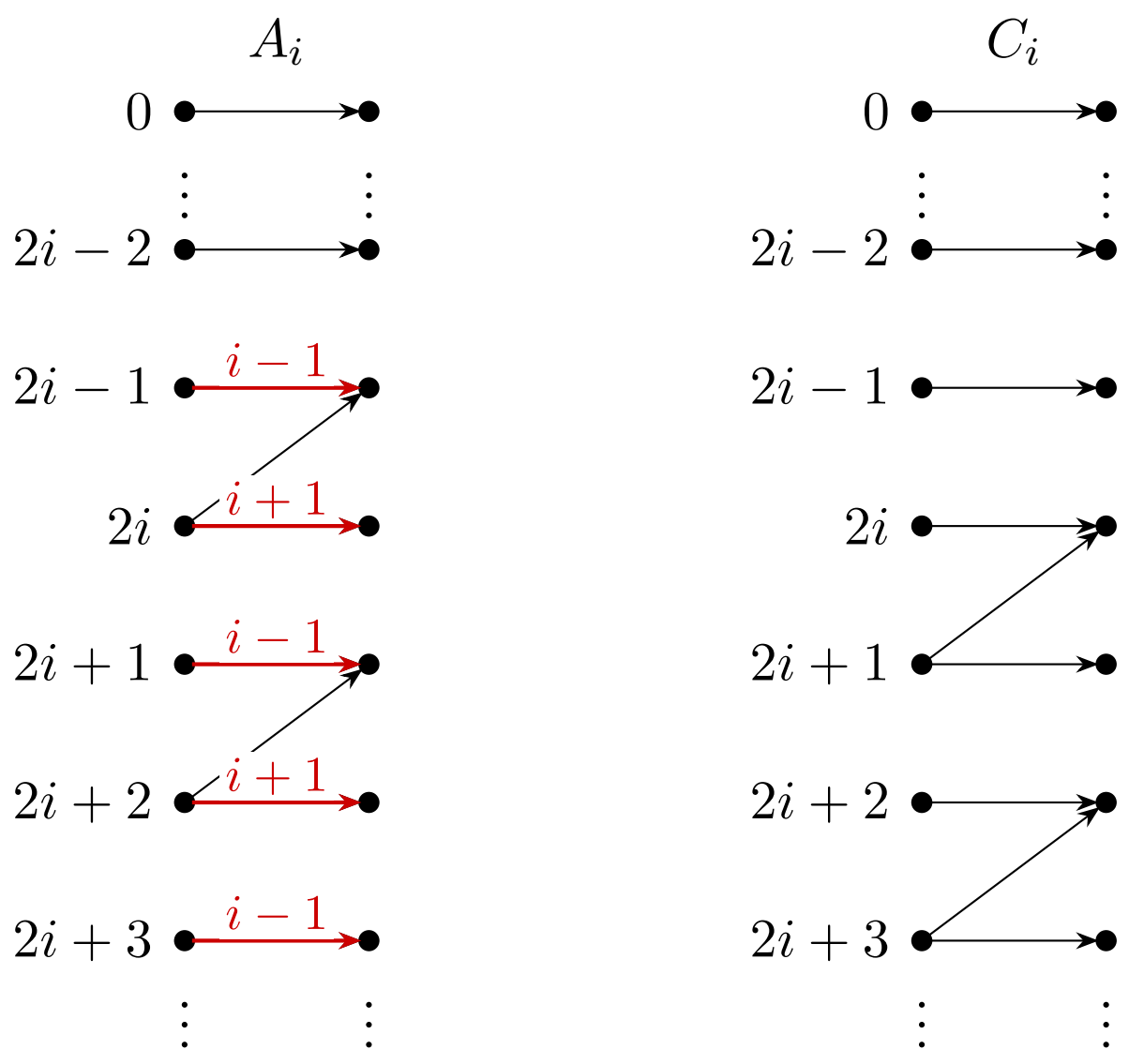

\begin{proof}
In Lemma \ref{lem:L}, we obtained $L=A_1C_1A_2C_2\cdots$.
It is easy to verify that the bidiagonal matrices $A_i$ and $C_i$ are the path matrices of the two digraphs shown in Figure \ref{fig:network-AC}, respectively.
Fix $N\ge1$.
Concatenating these digraphs in the order $A_1,C_1,A_2,C_2,\dots,A_N,C_N$ gives a digraph $\cD_N$ with the path matrix
$$
M_N=A_1C_1A_2C_2\cdots A_NC_N.
$$
In the network truncated at $x=2N$, extend each sink $V_s=(s,s)$, $0\le s\le2N$, horizontally to $W_s=(2N,s)$ by edges of weight $1$.
A path ending at $V_s$ extends uniquely along these auxiliary edges, with its weight unchanged.
So the extended digraph also has path matrix $M_N$.
Moreover, all factors of $L=A_1C_1A_2C_2\cdots A_NC_N\cdots$ after $C_N$ fix columns $0,\dots,2N$.
Hence, $(M_N)_{r,s}=L_{r,s}$ for $s\le2N$.
Therefore, taking the even-indexed sources and sinks, we derive
\[
P_{\cD}(U_{2n}\to V_{2k})=(M_N)_{2n,2k}=L_{2n,2k}=Q_{n,k}
\]
for $n\ge0$ and $0\le k\le N$.
Since $N$ is arbitrary, the result follows.
\end{proof}

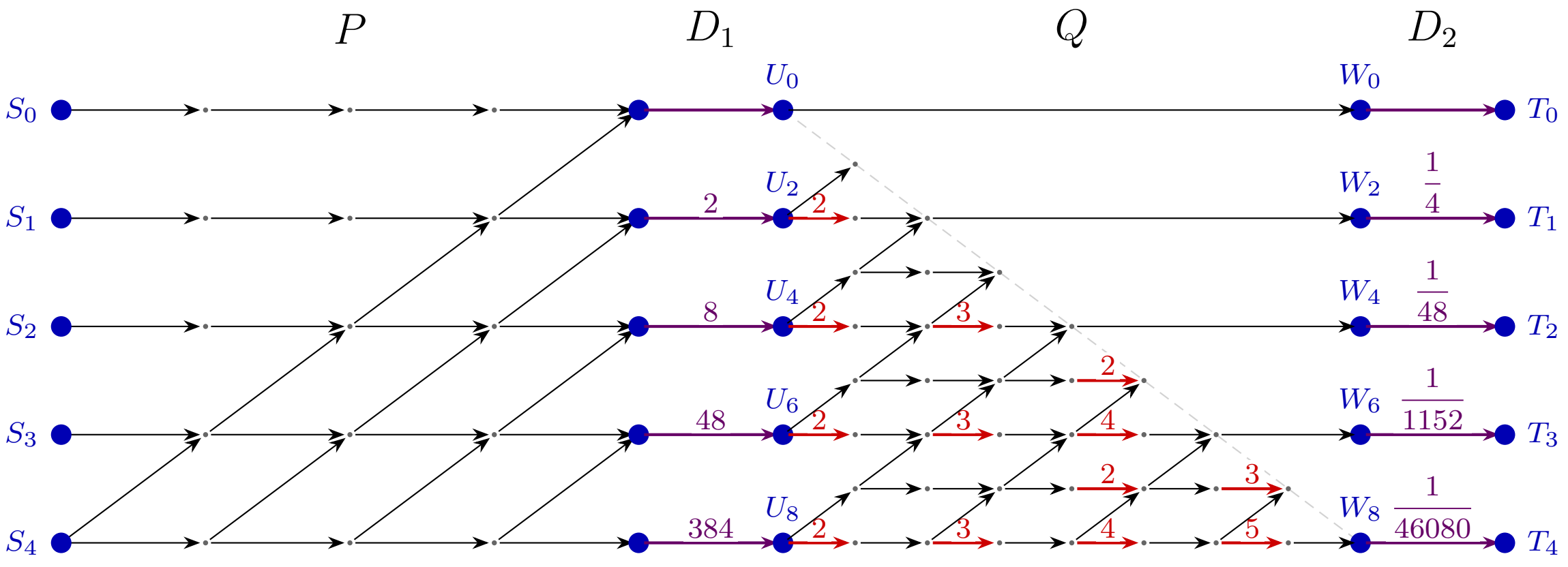
\begin{figure}[htbp]
	\centering
	\begin{tikzpicture}[x=1cm,y=-0.75cm,scale=0.7,
		vertex/.style={circle,fill=black,inner sep=1.5pt},,
		edge/.style={-{Stealth[length=1.6mm]},shorten <=1.5pt,shorten >=1.5pt},
		nonunit/.style={edge,red!80!black,thick},
		diagonal factor/.style={edge,violet!80!black,thick},
		weight/.style={fill=white,inner sep=1.1pt,font=\scriptsize},
		port label/.style={font=\scriptsize,inner sep=1pt}
		]
		\begin{scope}[x=2cm]
			\node at (2,-1.5) {$P$};
			\foreach \y in {0,2,...,8}{
				\foreach \x in {0,...,4}{
					\fill[black!60] (\x,\y) circle[radius=1pt];
				}
			}
			\foreach \x in {0,...,3}{
				\pgfmathtruncatemacro{\xp}{\x+1}
				\pgfmathtruncatemacro{\threshold}{4-\x}
				\foreach \r in {0,...,4}{
					\pgfmathtruncatemacro{\y}{2*\r}
					\draw[edge] (\x,\y)--(\xp,\y);
					\ifnum\r<\threshold\relax
						\else
						\pgfmathtruncatemacro{\yp}{2*(\r-1)}
						\draw[edge] (\x,\y)--(\xp,\yp);
					\fi
				}
			}
			\foreach \r in {0,...,4}{
				\def\sourcecolor{blue!70!black}
				\node[circle,fill=\sourcecolor,inner sep=2pt] at (0,2*\r) {};
				\node[port label,left=5pt,text=\sourcecolor] at (0,2*\r) {$S_{\r}$};
				\node[circle,fill=\sourcecolor,inner sep=2pt] at (4,2*\r) {};
			}
		\end{scope}
		\begin{scope}[xshift=8cm]
			\node at (1,-1.5) {$D_1$};
			\foreach \r/\w in {0/1,1/2,2/8,3/48,4/384}{
				\ifnum\r=0\relax
					\draw[diagonal factor] (0,0)--(2,0);
					\else
					\draw[diagonal factor] (0,{2*\r})--node[weight,above] {$\w$}(2,{2*\r});
				\fi
			}
		\end{scope}
		\begin{scope}[xshift=10cm]
			\node at (4,-1.5) {$Q$};
			\draw[gray!35,densely dashed] (0,0)--(8,8);
			\foreach \y in {1,...,8}{
				\foreach \x in {1,...,\y}{
					\fill[black!60] (\x,\y) circle[radius=1pt];
				}
			}
			\foreach \r in {0,...,8}{
				\pgfmathtruncatemacro{\rparity}{mod(\r,2)}
				\ifnum\rparity=0\relax
					\def\sourcecolor{blue!70!black}
					\def\sinkcolor{blue!70!black}
					\node[circle,fill=\sourcecolor,inner sep=2pt] at (0,\r) {};
					\node[port label,above=5pt,text=\sourcecolor] at (0,\r) {$U_{\r}$};
					\node[circle,fill=\sourcecolor,inner sep=2pt] at (8,\r) {};
					\node[port label,above=5pt,text=\sourcecolor] at (8,\r) {$W_{\r}$};
				\fi
			}
			\foreach \x in {0,...,7}{
				\foreach \y in {0,...,8}{
					\ifnum\y>\x\relax
						\pgfmathtruncatemacro{\xp}{\x+1}
						\pgfmathtruncatemacro{\xparity}{mod(\x,2)}
						\pgfmathtruncatemacro{\yparity}{mod(\y,2)}
						\ifnum\xparity=1\relax
							\def\edgeweight{1}
						\else
							\ifnum\yparity=0\relax
								\pgfmathtruncatemacro{\edgeweight}{\x/2+2}
							\else
								\pgfmathtruncatemacro{\edgeweight}{\x/2}
							\fi
						\fi
						\ifnum\edgeweight>0\relax
							\ifnum\edgeweight=1\relax
								\draw[edge] (\x,\y)--(\xp,\y);
							\else
								\draw[nonunit] (\x,\y)--node[weight,above] {$\edgeweight$}(\xp,\y);
							\fi
						\fi
						\ifnum\xparity=\yparity\relax
							\pgfmathtruncatemacro{\ym}{\y-1}
							\ifnum\ym<\xp\relax
							\else
								\draw[edge] (\x,\y)--(\xp,\ym);
							\fi
						\fi
					\fi
				}
			}
			\foreach \r in {0,...,3}{
				\draw[edge] (2*\r,2*\r)--(8,2*\r);	
			}
		\end{scope}
		\begin{scope}[xshift=18cm]
			\node at (1,-1.5) {$D_2$};
			\foreach \k/\den in {0/1,1/4,2/48,3/1152,4/46080}{
				\ifnum\k=0\relax
					\draw[diagonal factor] (0,0)--(2,0);
				\else
					\draw[diagonal factor] (0,{2*\k})--node[weight,above] {$\dfrac{1}{\den}$}(2,{2*\k});
				\fi
			}
			\foreach \r in {0,...,4}{
				\def\sourcecolor{blue!70!black}
				\node[circle,fill=\sourcecolor,inner sep=2pt] at (2,2*\r) {};
				\node[port label,right=5pt,text=\sourcecolor] at (2,2*\r) {$T_{\r}$};
			}
		\end{scope}
	\end{tikzpicture}
	\caption{A digraph for the Boolean-Narayana triangle restricted to $0\le n,k\le 4$.}
	\label{fig:BN-network}
\end{figure}

\begin{rem}
	Combining Lemma \ref{lem:BN-factor} and Corollary \ref{cor:Q-network}, we obtain a digraph for the Boolean-Narayana triangle $\BN=PD_1QD_2$ (see Figure \ref{fig:BN-network} for example),
	in which the four segments represent the digraphs for the Pascal triangle $P$, diagonal matrix $D_1$, matrix $Q$, and diagonal matrix $D_2$, respectively.
\end{rem}

\section{Gamma-positivity and strong $q$-log-convexity}

A real polynomial $f(q)$ is called \defn{$\ga$-positive} if it can be written as
\[
f(q)=\sum_{i=0}^{\floor{n/2}} \ga_i q^i (1+q)^{n-2i}
\]
with $\ga_i\ge0$ for all $i$.
The scalars $\ga_0,\ga_1,\ga_2,\dots$ are called the \defn{$\ga$-coefficients} of $f(q)$.
Gamma-positive polynomials arise often in combinatorics and geometry;
see the survey articles \cite{Bra15,Ath18}.

It is known \cite[p. 82]{Pet15} that if a polynomial with  symmetric nonnegative coefficients has only real zeros,
then it is $\ga$-positive.
B\'ona \cite{Bon26a} showed that the Boolean-Narayana polynomials $\BN_n(q)$ have symmetric nonnegative coefficients and only real zeros.
Thus, $\BN_n(q)$ is $\ga$-positive.
We recover this result by providing a combinatorial interpretation of the $\ga$-coefficients.

\begin{thm}\label{thm:BN-gpos}
	For $n\ge0$, the Boolean-Narayana polynomials $\BN_n(q)$ are $\ga$-positive with the expansion
	\begin{equation}\label{eq:BN-gaexp}
		\BN_n(q)=\sum_{j=0}^{\floor{n/2}} 2^j C_j \binom{n}{2j} q^j (1+q)^{n-2j}.
	\end{equation}
	Moreover, the $\ga$-coefficient $2^j C_j \binom{n}{2j}$ counts the 0-1-trees on $n+1$ vertices with exactly $j$ vertices having two children, in which every unique child is a left child.
\end{thm}

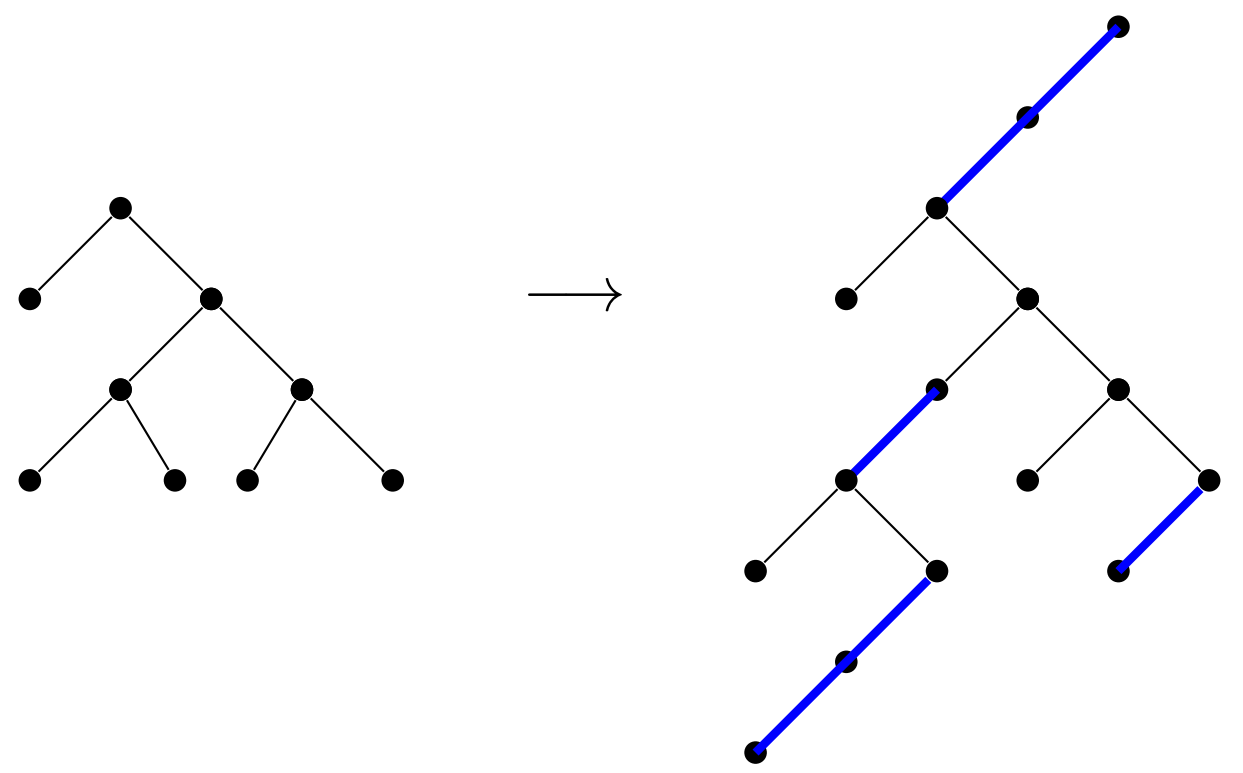
\begin{figure}[htbp]
	\centering
	\begin{tikzpicture}[y=-1cm,scale=0.6,
		ppath/.style={line width=0.8pt,line cap=round,line join=round},
		vtx/.style={circle,fill=black,inner sep=1.5pt},
		blpath/.style={color=blue,ultra thick}
		]
		\begin{scope}
			\node[vtx] (a) at (0,0) {};
			\node[vtx] (b) at (-1,1) {};
			\node[vtx] (c) at (1,1) {};
			\draw (b)--(a)--(c);
		\end{scope}
		\begin{scope}[xshift=1cm,yshift=-1cm]
			\node[vtx] (a) at (0,0) {};
			\node[vtx] (b) at (-1,1) {};
			\node[vtx] (c) at (1,1) {};
			\draw (b)--(a)--(c);
		\end{scope}
		\begin{scope}[yshift=-2cm]
			\node[vtx] (a) at (0,0) {};
			\node[vtx] (b) at (-1,1) {};
			\node[vtx] (c) at (0.6,1) {};
			\draw (b)--(a)--(c);
		\end{scope}
		\begin{scope}[xshift=2cm,yshift=-2cm]
			\node[vtx] (a) at (0,0) {};
			\node[vtx] (b) at (-0.6,1) {};
			\node[vtx] (c) at (1,1) {};
			\draw (b)--(a)--(c);
		\end{scope}
		\node at (5,1) {$\longrightarrow$};
		\begin{scope}[xshift=9cm]
			\node[vtx] at (2,-2) {};
			\node[vtx] at (1,-1) {};
			\draw[blpath] (2,-2)--(1,-1)--(0,0);
			\node[vtx] (a) at (0,0) {};
			\node[vtx] (b) at (-1,1) {};
			\node[vtx] (c) at (1,1) {};
			\draw (b)--(a)--(c);
		\end{scope}
		\begin{scope}[xshift=10cm,yshift=-1cm]
			\node[vtx] (a) at (0,0) {};
			\node[vtx] (b) at (-1,1) {};
			\node[vtx] (c) at (1,1) {};
			\draw (b)--(a)--(c);
		\end{scope}
		\begin{scope}[xshift=8cm,yshift=-3cm]
			\draw[blpath] (0,0)--(1,-1);
			\node[vtx] (a) at (0,0) {};
			\node[vtx] (b) at (-1,1) {};
			\node[vtx] (c) at (1,1) {};
			\draw (b)--(a)--(c);
			\node[vtx] at (0,2) {};
			\node[vtx] at (-1,3) {};
			\draw[blpath] (c)--(0,2)--(-1,3);
		\end{scope}
		\begin{scope}[xshift=11cm,yshift=-2cm]
			\node[vtx] (a) at (0,0) {};
			\node[vtx] (b) at (-1,1) {};
			\node[vtx] (c) at (1,1) {};
			\draw (b)--(a)--(c);
			\node[vtx] at (0,2) {};
			\draw[blpath] (c)--(0,2);
		\end{scope}
	\end{tikzpicture}
	\caption{Inserting chains of left edges (in blue) into a complete binary tree.}
	\label{fig:insert}
\end{figure}

\begin{proof}
Rewriting \eqref{eq:BN-formula-2} we have
\[
\BN(n,k)=\sum_{j=0}^{\min\{k,n-k\}} \binom{n}{2j} \binom{n-2j}{k-j} 2^j C_j.
\]
Then summing over $k$ gives
\begin{align*}
	\BN_n(q)
	&=\sum_{k=0}^n \BN(n,k) q^k
	=\sum_{j=0}^{\floor{n/2}} 2^j C_j \binom{n}{2j} \sum_{k\ge j} \binom{n-2j}{k-j} q^k\\
	&=\sum_{j=0}^{\floor{n/2}} 2^j C_j \binom{n}{2j} q^j (1+q)^{n-2j}.
\end{align*}

Recall that $\BN(n,k)$ counts the number of 0-1-trees on $n+1$ vertices with $k$ right edges.
Given such a 0-1-tree $T$, we assign weight $1$ to its left edges and weight $q$ to its right edges.
Define the weight $w(T)$ as the product of the weights of all its edges.
Then we have
\begin{equation}\label{eq:wT-BN}
	\sum_T w(T)=\sum_{k=0}^n \BN(n,k) q^k=\BN_n(q),
\end{equation}
where the first summation ranges over all weighted 0-1 trees on $n+1$ vertices.

On the other hand, fix the number $j$ of vertices having two children.
Such a 0-1-tree has $j+1$ leaves and $(n+1)-j-(j+1)=n-2j$ \defn{unary vertices} (vertices with a unique child).
Independently switching the child of any unary vertex between left and right partitions the trees into equivalence classes.
For example, the first, second, third, and sixth trees in Figure \ref{fig:0-1-tree} are in the same class.
Each class has a unique representative in which the child of every unary vertex is a left child.
Each of the $j$ vertices with two children has one left edge (with weight $1$) and one right edge (with weight $q$).
Hence, they contribute $q^j$ to the weight.
Each of the $n-2j$ unary vertices has either a left edge or a right edge.
Thus, they contribute $(1+q)^{n-2j}$ to the weight.
It remains to count the representatives.

Let $T$ be a 0-1-tree on $n+1$ vertices with exactly $j$ vertices having two children, in which every unique child is a left child.
After deleting all the $n-2j$ unique left edges,
contracting their endpoints, and ignoring the binary labels,
$T$ reduces to a complete binary tree
(that is, binary plane tree in which every non-leaf vertex has exactly two children)
on $2j+1$ vertices.
It is well known that the number of such complete binary trees is the Catalan number $C_j$.
We then put these $n-2j$ left edges back into the $2j+1$ positions according to the following rules:
if the position is a leaf, place a chain of left edges below it;
otherwise, place a chain of left edges above it.
Each chain may have length zero.
(See Figure \ref{fig:insert} for example.)
Hence, the number of ways to place the edges equals the number of ways to choose $n-2j$ out of $2j+1$ positions, with repetition allowed,
which is $\binom{(2j+1)+(n-2j)-1}{n-2j}=\binom{n}{2j}$.
There are $2^j$ choices for the binary labels of vertices having two children.
Hence, summing the total weights gives
\begin{equation}\label{eq:wT-gamma}
	\sum_T w(T)=\sum_{j=0}^{\floor{n/2}} 2^j C_j \binom{n}{2j} q^j (1+q)^{n-2j}.
\end{equation}
Combining \eqref{eq:wT-BN} with \eqref{eq:wT-gamma} yields \eqref{eq:BN-gaexp},
which completes the proof.
\end{proof}

Next, we will prove that the sequence of Boolean-Narayana polynomials is strongly $q$-log-convex by identifying it with the zeroth column of a recursive matrix.

For two real polynomials $f(q)$ and $g(q)$,
write $f(q)\ge_q g(q)$ if the coefficients of $f(q)-g(q)$ are all nonnegative.
A polynomial sequence $(f_n(q))_{n\ge0}$ is called \defn{$q$-log-convex} if $f_{n-1}(q) f_{n+1}(q)\ge_q f_n(q)^2$ for all $n\ge1$.
It is called \defn{strongly $q$-log-convex} if
\[
f_{m-1}(q) f_{n+1}(q)\ge_q f_m(q) f_n(q)
\]
for all $n\ge m\ge1$.
Strong $q$-log-convexity implies $q$-log-convexity,
but the converse can fail.
Such polynomials have been studied extensively, see \cite{LW07,CTW10,CWY10,CWY11,Zhu13,LZ15,Zhu18,Zhu24} for instance.

Let $\sig=(s_k(q))_{k\ge0}$ and $\tau=(t_k(q))_{k\ge1}$ be two sequences of polynomials with nonnegative coefficients.
Define an infinite lower triangular matrix $[a_{n,k}(q)]_{n,k\ge0}$ by the recurrence
\begin{equation}\label{eq:rec-st}
	a_{n,k}(q)=a_{n-1,k-1}(q)+s_k(q) a_{n-1,k}(q)+t_{k+1}(q) a_{n-1,k+1}(q), \quad n\ge1,
\end{equation}
where $a_{0,0}(q)=1$ and $a_{n,k}(q)=0$ unless $0\le k\le n$.
Following Aigner \cite{Aig01}, we call $[a_{n,k}(q)]_{n,k\ge0}$ the \defn{recursive matrix} corresponding to $\sig,\tau$.

\begin{thm}
	The sequence of Boolean-Narayana polynomials $(\BN_n(q))_{n\ge0}$ is strongly $q$-log-convex.
\end{thm}

\begin{proof}
Let
\[
J(q)=
\begin{bmatrix}
	1 \\
	1+q & 1 \\
	2q  & 1+q & 1 \\
	    & 2q  & 1+q    & 1  \\
		&     & \ddots & \ddots & \ddots
\end{bmatrix}
\]
be a lower triangular matrix with diagonal entries $1$ and first and second subdiagonals $1+q$ and $2q$, respectively.
Define a matrix $B(q)=[B_{n,k}(q)]_{n,k\ge0}$ recursively by
\[
B(q)=
\begin{bmatrix}
	1 \\
	  & B(q)
\end{bmatrix}
J(q).
\]
It is easy to verify that $B(q)$ is a recursive matrix with $s_k(q)=1+q$ and $t_{k+1}(q)=2q$ for $k\ge0$,
or equivalently, $B(q)$ satisfies the recurrence
\begin{equation}\label{eq:Bnk-rec}
	B_{n,k}(q)=B_{n-1,k-1}(q)+(1+q)B_{n-1,k}(q)+2q B_{n-1,k+1}(q)
\end{equation}
with $B_{0,0}(q):=1$ and $B_{n,k}(q)=0$ unless $0\le k\le n$.

We claim that the zeroth column of $B(q)$ is precisely the sequence of Boolean-Narayana polynomials $(\BN_n(q))_{n\ge0}$.
Indeed, let $G_k(q,t)=\sum_{n\ge0}B_{n,k}(q)t^n$.
Then \eqref{eq:Bnk-rec} implies
\begin{align*}
	G_0(q,t)&=1+(1+q)tG_0(q,t)+2qtG_1(q,t),\\
	G_k(q,t)&=tG_{k-1}(q,t)+(1+q)tG_k(q,t)+2qtG_{k+1}(q,t), \quad k\ge1.
\end{align*}
Rearranging these equations gives
\begin{equation}\label{eq:G0}
	G_0(q,t)=\frac{1}{1-(1+q)t-2qt G_1(q,t)/G_0(q,t)},
\end{equation}
and 
\begin{equation}\label{eq:Gk}
	\frac{G_k(q,t)}{G_{k-1}(q,t)}=\frac{t}{1-(1+q)t-2qt G_{k+1}(q,t)/G_k(q,t)}, \quad k\ge1.
\end{equation}
Iteration of \eqref{eq:Gk} and \eqref{eq:G0} leads to the continued fraction
\[
G_0(q,t)=
\cfrac{1}{
	1-(1+q)t-
	\cfrac{2qt^2}{
		1-(1+q)t-
		\cfrac{2qt^2}{
			1-(1+q)t-\cdots}}}.
\]
Then we have
\[
G_0(q,t)=\frac{1}{1-(1+q)t-2qt^2G_0(q,t)},
\]
and hence,
\[
2qt^2G_0(q,t)^2-(1-(1+q)t)G_0(q,t)+1=0.
\]
A routine calculation gives
\begin{align}
	\nonumber
	G_0(q,t)
	&=\frac{1-(1+q)t-\sqrt{(1-(1+q)t)^2-8qt^2}}{4qt^2}\\
	\label{eq:G0qt}
	&=\frac{1}{1-(1+q)t} C\left(\frac{2qt^2}{(1-(1+q)t)^2}\right),
\end{align}
where $C(t)=\frac{1-\sqrt{1-4t}}{2t}$ is the generating function of the Catalan numbers.

On the other hand, let $F(q,t)=\sum_{n\ge0}\BN_n(q)t^n$.
Then by \eqref{eq:BN-gaexp}, we obtain
\begin{align}
	\nonumber
	F(q,t)
	&=\sum_{n\ge0} \sum_{j=0}^{\floor{n/2}} 2^j C_j \binom{n}{2j} q^j (1+q)^{n-2j} t^n\\
	\nonumber
	&=\sum_{j\ge0} C_j (2qt^2)^j \sum_{m\ge0} \binom{m+2j}{2j} ((1+q)t)^m\\
	\nonumber
	&=\sum_{j\ge0} \frac{C_j (2qt^2)^j}{(1-(1+q)t)^{2j+1}}\\
	\label{eq:Fqt}
	&=\frac{1}{1-(1+q)t} C\left(\frac{2qt^2}{(1-(1+q)t)^2}\right).
\end{align}
Comparing \eqref{eq:G0qt} with \eqref{eq:Fqt}, we have $G_0(q,t)=F(q,t)$.
Taking the coefficients of $t^n$ yields $B_{n,0}(q)=\BN_n(q)$, as desired.

Zhu \cite[Theorem 2.1]{Zhu13} showed that if $s_{k-1}(q) s_k(q)\ge_q t_k(q)$ for $k\ge1$ then the zeroth column of the recursive matrix
defined by \eqref{eq:rec-st} is strongly $q$-log-convex.
For the matrix $B(q)$ we have
\[
s_{k-1}(q) s_k(q)=(1+q)^2\ge_q 2q= t_k(q).
\]
Therefore, $(\BN_n(q))_{n\ge0}$, as the zeroth column of $B(q)$, is strongly $q$-log-convex.
\end{proof}



\section{Asymptotic normality}

Let $[a_{n,k}]_{n,k\ge0}$ be a lower triangular matrix with nonnegative entries and positive row sums.
Let
\[
p(n,k)=\frac{a_{n,k}}{\sum_{j=0}^n a_{n,j}}
\]
denote the normalized probabilities.
Following Bender \cite{Ben73},
we say that the numbers $a_{n,k}$ are \defn{asymptotically normal by a central limit theorem} if
\begin{equation}\label{eq:CLT}
    \lim_{n\to\infty} \sup_{x\in\R}
	\left|\sum_{k\le\mu_n+x\sigma_n} p(n,k)-\frac{1}{\sqrt{2\pi}} \int_{-\infty}^x e^{-t^2/2}dt\right|=0,
\end{equation}
where $\mu_n$ and $\sigma^2_n$ are the mean and variance of $(a_{n,k})$, respectively.
We say that $a_{n,k}$ are \defn{asymptotically normal by a local limit theorem} on $\R$ if
\begin{equation}\label{eq:LLT}
	\lim_{n\to\infty} \sup_{x\in\R}
	\left|\sigma_n p(n,\floor{\mu_n+x\sigma_n})-\frac{1}{\sqrt{2\pi}} e^{-x^2/2}\right|=0.
\end{equation}
In this case,
$$
a_{n,k}\sim\frac{e^{-x^2/2} \sum_{j=0}^n a_{n,j}}{\sigma_n \sqrt{2\pi}}
\text{ as } n\to\infty,
$$
where $k=\mu_n+x\sigma_n$ and $x=O(1)$.
Clearly, the validity of \eqref{eq:LLT} implies that of \eqref{eq:CLT}.

Many well-known combinatorial sequences satisfy central and local limit theorems.
For instance, the famous de Moivre-Laplace theorem states that the binomial coefficients are asymptotically normal (by central and local limit theorems).
Other examples include the Stirling numbers of both kinds \cite{Fel45,Har67},
the Eulerian numbers \cite{CKS72} with a generalization \cite{HCD20},
the Narayana numbers \cite{CMW20} with a generalization \cite{CYZ22},
and variations of the Motzkin numbers \cite{CWZ20}.
See Canfield \cite{Can15} for an excellent survey on such a topic.

A standard approach to demonstrating asymptotic normality is the following criterion,
which was first used by Harper \cite{Har67} to show the asymptotic normality of the Stirling numbers of the second kind
(see \cite[Theorem 2]{Ben73} for instance and \cite[Example 3.4.2]{Can15} for historical remarks).

\begin{lem}[{\cite[Theorem 2]{Ben73}}]\label{lem:AN}
	Suppose that $A_n(x)=\sum_{k=0}^n a_{n,k} x^k$ have only real roots,
	with $a_{n,k}\ge0$ and $A_n(1)>0$.
	Let
	\begin{equation}\label{eq:mean-var}
		\mu_n=\frac{A'_n(1)}{A_n(1)}, \quad
		\sigma^2_n=\frac{A''_n(1)}{A_n(1)}+\mu_n-\mu^2_n.
	\end{equation}
	If $\sigma^2_n\to+\infty$,
	then the numbers $a_{n,k}$ are asymptotically normal (by both central and local limit theorems)
	with mean $\mu_n$ and variance $\sigma_n^2$.
\end{lem}

\begin{thm}
	The Boolean-Narayana numbers $\BN(n,k)$ are asymptotically normal by both central and local limit theorems.
\end{thm}

\begin{proof}
For brevity, denote $B_n(q):=\BN_n(q)$.
It is known \cite[Theorem 5.1]{Bon26a} that the polynomials $B_n(q)$ have only real zeros.
Hence, by Lemma \ref{lem:AN}, it suffices to show that the variances $\sig_n^2$ of $\BN(n,k)$ tend to infinity.

Recall that
\[
B_n(q)=\sum_{j=0}^{\floor{n/2}} 2^j C_j \binom{n}{2j} q^j (1+q)^{n-2j}.
\]
Then we have
\[
B_n'(q)=\sum_{j=0}^{\floor{n/2}} 2^j C_j \binom{n}{2j} q^{j-1} (1+q)^{n-2j-1} [j+(n-j)q]
\]
and
\[
B_n''(q)=\sum_{j=0}^{\floor{n/2}} 2^j C_j\binom{n}{2j}q^{j-2}(1+q)^{n-2j-2} \left[j(j-1)+2j(n-j-1)q+(n-j)(n-j-1)q^2\right].
\]
Hence,
\[
B_n'(1)=n\sum_{j=0}^{\floor{n/2}} 2^{n-j-1} C_j \binom{n}{2j},
\]
and
\begin{align*}
	B_n''(1)
	&=\sum_{j=0}^{\floor{n/2}} 2^{n-j-2} C_j \binom{n}{2j} (n^2-n-2j)\\
	&=n(n-2)\sum_{j=0}^{\floor{n/2}} 2^{n-j-2} C_j \binom{n}{2j}+\sum_{j=0}^{\floor{n/2}} 2^{n-j-2} C_j \binom{n}{2j} (n-2j)\\
	&=n(n-2)\sum_{j=0}^{\floor{n/2}} 2^{n-j-2} C_j \binom{n}{2j}+n\sum_{j=0}^{\floor{n/2}} 2^{n-j-2} C_j \binom{n-1}{2j}.
\end{align*}
Note that
\[
B_n(1)=\sum_{j=0}^{\floor{n/2}} 2^{n-j} C_j \binom{n}{2j}.
\]
Thus, by \eqref{eq:mean-var}, we derive
\[
\mu_n=\frac{B'_n(1)}{B_n(1)}=\frac{n}{2}
\]
and
\[
\sigma^2_n=\frac{B''_n(1)}{B_n(1)}+\mu_n-\mu^2_n
=\frac{n(n-2)}{4}+\frac{n}{2}\cdot\frac{B_{n-1}(1)}{B_n(1)}+\frac{n}{2}-\frac{n^2}{4}
=\frac{n}{2}\cdot\frac{B_{n-1}(1)}{B_n(1)}.
\]
B\'ona \cite[Theorem 6.1]{Bon26a} gave the recurrence
\[
(n+2)B_n(q)=(2n+1)(1+q)B_{n-1}(q)-(n-1)(q^2-6q+1)B_{n-2}(q), \quad n\ge2,
\]
with $B_0(q)=1$ and $B_1(q)=1+q$.
Plugging in $q=1$ and dividing both sides by $(n+2)B_{n-1}(1)$ yields
\begin{align*}
	\frac{B_n(1)}{B_{n-1}(1)}
	&=\frac{2(2n+1)}{n+2}+\frac{4(n-1)}{n+2}\cdot\frac{B_{n-2}(1)}{B_{n-1}(1)}\\
	&=2+\frac{2(n-1)}{n+2}\left(1+\frac{2B_{n-2}(1)}{B_{n-1}(1)}\right).
\end{align*}
It is easy to see that $B_n(1)/B_{n-1}(1)\ge2$ for $n\ge1$, since every $B_n(1)$ is positive.
Hence, we obtain
\[
\frac{B_n(1)}{B_{n-1}(1)}
\le \frac{2(2n+1)}{n+2}+\frac{4(n-1)}{n+2}\cdot\frac{1}{2}
=\frac{6n}{n+2}.
\]
It follows that
\[
\sigma^2_n=\frac{n}{2}\cdot\frac{B_{n-1}(1)}{B_n(1)}
\ge\frac{n+2}{12}
\to+\infty.
\]
Therefore, by Lemma \ref{lem:AN}, the Boolean-Narayana numbers are asymptotically normal.
\end{proof}

\section{An alternative combinatorial definition}

We construct a bijection between 0-1-trees and colored Motzkin paths,
hence present an alternative combinatorial definition of the Boolean-Narayana numbers.

A \defn{Motzkin path of length $n$} is a lattice path from $(0,0)$ to $(n,0)$ that never goes below the horizontal axis and
consists of up steps $U=(1,1)$, horizontal steps $H=(1,0)$, and down steps $D=(1,-1)$.
Following \cite{AM06}, a \defn{colored Motzkin path} is a Motzkin path in which each up step has color $0$ or $1$, each horizontal step has color $L$ or $R$, and down steps are uncolored.

\begin{thm}\label{thm:tree-Motzkin}
	For every $n\ge 0$, there is a bijection between 0-1-trees on $n+1$ vertices and colored Motzkin paths of length $n$.
\end{thm}

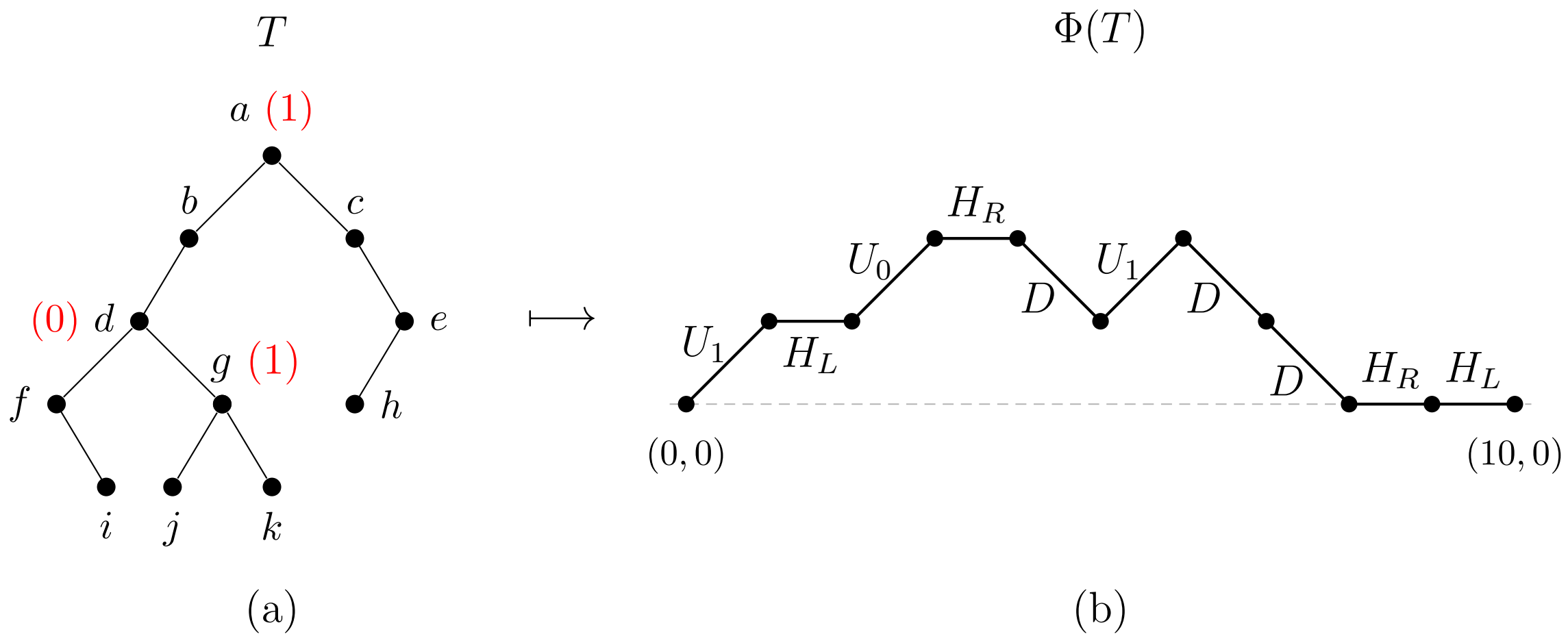
\begin{figure}[htbp]
	\centering
	\begin{tikzpicture}[scale=0.8,
		vtx/.style={circle,fill=black,inner sep=1.8pt},
		Mstep/.style={line width=0.8pt,line cap=round,line join=round},
		every label/.style={font=\small}
		]
		\begin{scope}
			\node[vtx,label=above:{$a~\tre{(1)}$}] (a) at (0,0) {};
			\node[vtx,label=above:{$b$}] (b) at (-1,-1) {};
			\node[vtx,label=above:{$c$}] (c) at (1,-1) {};
			\draw (b)--(a)--(c);
		\end{scope}
		
		\begin{scope}[xshift=-1.6cm,yshift=-2cm]
			\node[vtx,label=left:{$f$}] (f) at (-1,-1) {};
			\node[vtx,label=left:{\tre{(0)}~$d$}] (d) at (0,0) {};
			\node[vtx] (c) at (1,-1) {};
			\draw (f)--(d)--(c);
		\end{scope}

		\begin{scope}[xshift=-0.6cm,yshift=-3cm]
			\node[vtx,label=below:{$j$}] (j) at (-.6,-1) {};
			\node[vtx] (g) at (0,0) {};
			\node at (0.4,0.5) {$g$~\tre{(1)}};
			\node[vtx,label=below:{$k$}] (k) at (.6,-1) {};
			\draw (j)--(g)--(k);
		\end{scope}
		\draw (-1,-1)--(-1.6,-2) (-2.6,-3)--(-2,-4) (1,-1)--(1.6,-2)--(1,-3);
		\node[vtx,label=right:{$e$}] at (1.6,-2) {};
		\node[vtx,label=right:{$h$}] at (1,-3) {};
		\node[vtx,label=below:{$i$}] at (-2,-4) {};
		\node at (0,1.5) {$T$};
		\node at (0,-5.5) {(a)};

		\node at (3.5,-2) {$\longmapsto$};

		\begin{scope}[xshift=5cm,yshift=-3cm]
			\draw[gray!55,densely dashed] (-.2,0)--(10.2,0);
			\foreach \t/\h in {0/0,1/1,2/1,3/2,4/2,5/1,6/2,7/1,8/0,9/0,10/0}
				\coordinate (P\t) at (\t,\h);
			\foreach \t in {0,...,10}
				\fill (P\t) circle(0.1);
			\draw[Mstep] (P0)--node[above left=-6pt] {$U_1$}(P1);
			\draw[Mstep] (P1)--node[below] {$H_L$}(P2);
			\draw[Mstep] (P2)--node[above left=-6pt] {$U_0$}(P3);
			\draw[Mstep] (P3)--node[above] {$H_R$}(P4);
			\draw[Mstep] (P4)--node[below left=-5pt] {$D$}(P5);
			\draw[Mstep] (P5)--node[above left=-6pt] {$U_1$}(P6);
			\draw[Mstep] (P6)--node[below left=-5pt] {$D$}(P7);
			\draw[Mstep] (P7)--node[below left=-5pt] {$D$}(P8);
			\draw[Mstep] (P8)--node[above] {$H_R$}(P9);
			\draw[Mstep] (P9)--node[above] {$H_L$}(P10);
			\node[below=5pt] at (P0) {\footnotesize $(0,0)$};
			\node[below=5pt] at (P10) {\footnotesize $(10,0)$};
			\node at (5,4.5) {$\Phi(T)$};
			\node at (5,-2.5) {(b)};
		\end{scope}
	\end{tikzpicture}
	\caption{(a) A 0-1 tree on $11$ vertices, where the vertices with two children have the binary labels shown in red parentheses.
	(b) The corresponding colored Motzkin path of (a).}
	\label{fig:tree-Motzkin}
\end{figure}

\begin{proof}
Let $T$ be a 0-1-tree on $n+1$ vertices.
Traverse the vertices of $T$ in preorder:
visit the root first, then recursively traverse the left subtree and the right subtree.
For instance, traversing the tree in Figure \ref{fig:tree-Motzkin}(a) in preorder gives $abdfigjkceh$.
For each vertex in $T$, record a colored Motzkin step according to the following rules:

\begin{center}
\begin{tabular}{c|c}
	Vertex type & Motzkin step \\
	\hline
	labeled $0$, two children & $U_0$ \\
	labeled $1$, two children & $U_1$ \\
	only a left child & $H_L$ \\
	only a right child & $H_R$ \\
	leaf & $D$
\end{tabular}
\end{center}
where the subscripts in the Motzkin steps indicate the color.
By definition, the last vertex in preorder must be a leaf.
Delete its corresponding down step and denote the resulting
path by $\Phi(T)$;
see Figure \ref{fig:tree-Motzkin} as an example.

We first show that $\Phi$ is well defined, i.e., $\Phi(T)$ is a colored Motzkin path of length $n$. 
Recall that in a 0-1-tree, each vertex with two children increases the number of leaves by one.
Hence, the number of leaves is exactly one more than the number of vertices with two children.
Then, after deleting the last leaf, we obtain a lattice path which has the same number of up and down steps.
So $\Phi(T)$ ends at $(n,0)$.
Moreover, before a leaf occurs, there must be a corresponding vertex having two children,
since the preorder visits parents before children.
Let $h_m$ be the height after the first $m$ ($1\le m\le n$) steps of $\Phi(T)$.
Then $h_m\ge0$, as in the first $m$ steps,
the number of up steps (vertices with two children) is no less than the number of down steps (leaves).
It follows that $\Phi(T)$ never goes below the horizontal axis.
Thus, $\Phi(T)$ is a Motzkin path.

Next, we show that $\Phi$ is a bijection.
Let $M$ be a colored Motzkin path of length $n$.
Append a down step to the end of $M$ and read its steps from left to right.
We always draw a 0-1-tree $T$ in preorder,
and define the \defn{pending list} of $T$ as an ordered list of vertices waiting to be visited.
Initially, the pending list contains only the root position.
For each of the $n+1$ Motzkin steps, we remove the first position from the pending list of $T$ and place a vertex there
according to the rules in the above table.
If the vertex is not a leaf, add the position(s) of its child(ren) to the front of the list;
if it has two children, put the left position before the right one.
Note that each up (resp. down) step of $M$ increases (resp. decreases) the number of pending positions by one, and each horizontal step does not change this number.
Hence, the number of pending positions always equals the current height of $M$ plus one.
This ensures that a position is available before every step,
and the pending list becomes empty after the appended down step.
Therefore, this procedure yields a unique 0-1-tree on $n+1$ vertices and inverts $\Phi$.
\end{proof}

\begin{cor}\label{cor:BN-Motzkin}
The Boolean-Narayana number $\BN(n,k)$ counts the colored Motzkin paths of length $n$
with exactly $k$ steps that are either up steps or horizontal steps of color $R$.
\end{cor}

\begin{proof}
By Definition \ref{defi:BN}, the numbers $\BN(n,k)$ count 0-1-trees on $n+1$ vertices with $k$ right edges.
Note that the right edges are contributed by both vertices with two children and vertices with one right child.
Under the bijection $\Phi$ constructed in Theorem~\ref{thm:tree-Motzkin},
each vertex with two children corresponds to an up step,
while a vertex with only a right child corresponds to a horizontal step with color $R$.
Thus the result follows.
\end{proof}

In terms of the colored Motzkin paths, we obtain another interpretation of the $\ga$-coefficients of the Boolean-Narayana polynomials $\BN_n(q)$.

\begin{cor}\label{cor:gamma-Motzkin}
	For $n\ge0$ and $0\le j\le\floor{n/2}$, 
	the $\ga$-coefficient $2^j C_j\binom{n}{2j}$ of $\BN_n(q)$ counts the colored Motzkin paths of length $n$ with exactly $j$ up steps, in which every horizontal step has color $L$.
\end{cor}

\begin{proof}
For a colored Motzkin path $M$, assign weight $q$ to its up steps and horizontal steps with color $R$, and weight $1$ to all other steps.
Define the weight of a path as the product of the weights of all its steps.
Then, by Corollary \ref{cor:BN-Motzkin},
\[
\sum_{M} w(M)=\sum_{k=0}^n \BN(n,k)q^k=\BN_n(q),
\]
where the first summation ranges over all colored Motzkin paths of length $n$.

We say that two colored Motzkin paths are in the same equivalence class if they are identical except for the colors of horizontal steps.
In each class, we choose the unique path whose horizontal
steps are all colored $L$ as its representative.
Suppose that this representative has $j$ up steps.
Then it has $j$ down steps and $n-2j$ horizontal steps.
Since each horizontal step may have color $L$ or $R$, the paths in this class contribute $q^j(1+q)^{n-2j}$ in total to the weight.
So the $\ga$-coefficient is the number of equivalence classes according to $j$.

It remains to count the representatives.
After deleting all the horizontal steps, each Motzkin path reduces to a Dyck path consisting of $j$ up steps and $j$ down steps.
There are $C_j$ such Dyck paths.
Then, we insert the $n-2j$ horizontal steps with color $L$ into the Dyck path.
There are $\binom{n}{n-2j}=\binom{n}{2j}$ choices for their positions.
Note that each up step has color $0$ or $1$, giving $2^j$ choices.
Therefore, the number of representatives with exactly $j$ up steps is $2^j C_j\binom{n}{2j}$.
This completes the proof.
\end{proof}


\end{document}